\documentclass[12pt]{amsart}

\usepackage{amssymb,latexsym, amsfonts}

\newtheorem{theorem}{Theorem}[section]

\title[Sandwich Theorem] {The Sandwich Theorem for Sublinear and Superlinear Functionals }

\author{A.T. Diab, S. I. Nada, D. L. Fearnley}

\address{Dept. of Math., Faculty of Science \\
Ain Shams University \\
Cairo \\ Egypt}

\address{Dept. of Math., Faculty of Science \\
Qatar University \\
Doha \\ Qatar}

\address{Dept. of Math. \\
Utah Valley University \\
Orem, Utah \\ USA}

\email{adeldiab80@hotmail.com}
\email{snada@qu.edu.qa}
\email{fearnlda@uvu.edu}


\begin{document}
\begin{abstract}

   The Hahn-Banach theorem is an extension theorem for linear
functionals which preserves certain properties. Specifically, if a
linear functional is defined on a subspace of a real vector space
which is dominated by a sublinear functional on the entire space,
then this functional can be extended to a linear functional on the
entire space which is still dominated by a sublinear functional. In
this paper, we generalize this result to show that a linear
functional defined on a subspace of a real vector space which is
dominated by a sublinear functional and also dominates a superlinear
functional on the entire space can be extended to a linear
functional on the entire space which is also dominated by a
sublinear functional and dominates a superlinear functional.

\end{abstract}

\maketitle \pagenumbering{arabic} \setcounter{page}{1}
\section{Introduction}

   Extensions of linear functionals have been studied in locally convex
spaces by several authors. In fact, the importance of the
Hahn-Banach theorem arises from its wide variety of applications,
including complex and functional analysis and thermodynamics. The
Hahn-Banach theorem has important implications for convex sets, and
is the foundation for an effective treatment of optimization. It is
also used to solve many problems in linear algebra, conic duality
theory, the minimax theorem, piecewise approximation of convex
functionals, extensions of positive linear functionals, and other
results from modern control [1,5,8]. The Hahn-Banach theorem was
proved by Hans Hahn (1879- 1934) in 1927 and later by Stefan Banach
(1892-1945) in 1929. It states that if $E$ is a real vector space,
$M$ is a subspace of $E$ and $S$ is a sublinear functional on $E$,
and $f_0$ is a linear functional on $M$  such that $ f_0(x)\leq
S(x)$ for every $ x\in M$ (or equivalently, $ -S(-x)\leq f_0(x)\leq
S(x)$), then there exists a linear functional $L$ on $E$ such that
$L|_M =f_0$ and $ L(x)\leq S(x)$ for every $ x\in E$ (or
equivalently, $ -S(-x)\leq L(x)\leq S(x)$). We will show that if $P$
is a superlinear functional on $E$ satisfying the condition that
$P(x)\leq f_0(x)\leq S(x)$ for every $ x\in M$ such that $f_0(x)\le
T(x)$ for every $x\in M$, where $T(x)=\inf\limits_{y\in
E}\{S(x+y)-P(y)\}$ for every $x\in E$, then there exists a linear
functional $L$ on $E$ such that $L|_M =f_0$ and $ P(x)\leq L(x)\leq
S(x)$ for every $ x\in E$. We will refer to this as the sandwich
theorem for a sublinear and superlinear functionals. Note that if we
take $ P(x)=-S(-x)$ for every $ x\in E$, then we obtain the
Hahn-Banach theorem. As an application of the Hahn-Banach theorem we
have the following theorem: There exists a linear functional $L$ on
$E$ for every sublinear functional $S$ on $E$ such that $L(x)\leq
S(x)$ for every $ x\in E$ (we refer to this as the classical
Hahn-Banach theorem). We will also use the sandwich theorem for a
sublinear and superlinear functionals to show that, for every a
sublinear functional $S$ on $E$ and for every superlinear functional
$P$ on $E$, if $P(x)\leq S(x)$ for every $ x\in E$ then there exists
a linear functional $L$ on $E$ such that $P(x)\leq L(x)\leq S(x)$
for every $ x\in E$. Note that if we take $P(x)=-S(-x)$ for every $
x\in E$, then this result implies the classical Hahn-Banach theorem.

\section{Relations between sublinear and superlinear functionals defined on a real linear space}

  In this section we present some basic definitions and results from functional analysis
concerning sublinear and superlinear functionals. In the following
definitions we suppose that $E$ is a real vector space [6,7,9].

\

\noindent \textbf{Definition 2.1.}  Let $E$ be a real vector space,
$M$ be a subspace of $E$ and $f_0$ be a linear functional on $M$.
Then a linear functional $L$ on $E$ is called an \textit{extension
functional} of $f_0$ if $L(x) = f_0(x)$ for every $x \in M$ (denoted
$L|_M =f_0$).

\

 \noindent \textbf{Definition 2.2.}  The functional $
S: E \rightarrow \mathbb{R}$ is called \textit{sublinear} if it possesses the
following properties:

\

 \noindent (1)  $S(x+y)\leq S(x)+S(y)$ for every $x,y\in E$ (i.e., $S$ is
 subadditive),\\
 (2)  $S(\alpha x)= \alpha S(x)$ for every $ x\in E$ and $\alpha >0$ (i.e., $S$ is positively homogeneous). \
 The set of all sublinear functionals on $E$ is denoted by $Subl(E)$.

 \

\noindent \textbf{Definition 2.3.}  The functional $ P: E
\rightarrow \mathbb{R}$ is called \textit{superlinear} if it possesses the
following properties:

\

\noindent (1)$P(x+y) \geq P(x)+P(y)$ for every $ x, y\in E$ (i.e., $P$ is
superadditive),\\
(2)$ P(\alpha x)= \alpha P(x)$ for every $ x\in E$ and $\alpha >0$
(i.e., $P$ is positively homogeneous).\ The set of all superlinear
functionals on $E$ is denoted by $Supl(E)$.

\

Note that $Lin(E)= Subl(E)\bigcap Supl(E)$, where $Lin(E)$ is the
set of all linear functionals on $E$. Also, it is easy to verify
that if $E\equiv\mathbb{R}^{n}, S(x)=\sqrt{\sum\limits_{i=1}^n
x_i^2}  = \|x\|$ and  $P(x)=-\sqrt{\sum\limits_{i=1}^{n-1}
x_i^2}+x_n$ for every $x\in \mathbb{R}^{n}$, then $ S\in Subl(E)$
and $P\in Supl(E)$.
Also, $S\in Subl(E)$ if and only if $-S\in Supl(E)$.

\

Throughout the this paper we will use the symbols $S$ and $P$ to
denote sublinear and superlinear functionals on the real vector
space $E$.

\

\noindent \textbf{Lemma 2.1.}  \textit{(i) $S( \underline{0})= 0$ and
$P(\underline{0}) = 0$, where $\underline{0}$ is the zero vector of
the real vector space $E$.\\ (ii) $-S(-x)\le S(x)$ and $-P(-x)\geq
P(x)$ for every $x\in E$.\\ (iii) A functional $f$ on $E$ is linear
if and only if it is additive and positively homogeneous.}

\

\noindent \textbf{Proof.} (i) Since sublinear functionals are
positively homogeneous, $S(\underline{0})=S(2 \underline 0) =
2.S(\underline 0+ \underline0)=S(\underline 0)+S(\underline 0)$.
Therefore $S(\underline 0)= 0$.  Similarly, $P(\underline{0})= 0$.\\
(ii) Since $0 =S(\underline{0})=S(-x+x)\le S(-x)+S(x)$ for every
$x\in E$, then $-S(-x)\le S(x)$. Also,
$0=P(\underline{0})=P(-x+x)\geq P(-x)+P(x)$ for every $x\in E$, so
$-P(-x)\geq P(x)$ for every $x\in E$.\\
(iii) By definition, linear functionals are additive and positively
homogeneous. Let $f$ be an additive positively homogeneous
functional on $E$.  It remains to show that $f(\alpha x)=\alpha
f(x)$ for every $x \in E$ and $\alpha<0$. Since $f$ is additive,
$f(\underline{0}+\underline{0})=f(\underline{0})+f(\underline{0})$.
Thus, $f(\underline{0})=0$, and consequently $f(-x)=-f(x)$ for every
$x\in E$. Let $x\in E$ and $\alpha =-t$, where $t>0$. Then $f(\alpha
x)= f(-t x)=t f(-x)=-t f(x)=\alpha f(x)$. Thus $f$ is homogeneous
and $f$ is linear.

\

\noindent \textbf{Lemma 2.2.}  \textit{If $S(x)\le P(x)$ for every $x\in E$,
then $P$ and $S$ are linear functionals on $E$. Moreover $S\equiv P$
on $E$.}

\

\noindent \textbf{Proof.}  Since $S$ is sublinear and $P$ is
superlinear, then from(ii)in Lemma 2.1, we obtain that $-S(-x)\le
S(x)$ and  $-P(-x)\ge P(x)$ for every $ x\in E$. So, $S(-x)\le
P(-x)\le -P(x)\le - S(x)$. Thus, $S(-x)=-S(x)$ for every $x\in E$.
For all $x,y\in E$, $S(x+y)=-S(-x-y)\ge -S(-x)-S(-y)$.  Hence, $S$
is both subadditive and superadditive and is therefore additive.
Since $S$ is positively homogeneous and additive, from (iii) in
Lemma 2.1 we can conclude that $S$ is linear. By a similar argument,
$P$ is linear. Finally, from (ii) in Lemma 2.1 we see that $-
S(x)\le S(-x)\le P(-x)\le -P(x)$. Hence, $P(x)\le S(x)$ for every
$x\in E$. Therefore $P(x) = S(x)$ for every $x\in E$.

\

\noindent \textbf{Lemma 2.3.}  \textit{If $f$ is a homogeneous functional on
$E$ and $g$ is an arbitrary functional on $E$, then the following
statements are equivalent: \newline (i)  $f(x)\le g(x)$ for every
$x\in E$, \newline (ii) $-g(-x)\le f(x)\le g(x)$ for every $x\in E$.}

\

\noindent \textbf{Proof.} Clearly, (ii) implies (i). Let $f(x)\le
g(x)$ for every $x\in E$. Then $f(-x)\le g(-x)$ and
consequently$-f(x)=f(-x)\le g(-x)$ for every $x\in E$. Therefore
$-g(-x)\le f(x)$ and hence $-g(-x)\le f(x)\le g(x)$.

\

\noindent \textbf{Corollary 2.1.}   \textit{If $f$ is a linear functional on
$E$ and  $f(x)\le S(x)$ for every $x\in E$, where $S$ is a sublinear
functional on $E$, then $-S(-x)\le f(x)\le S(x)$ for every $x\in E$.}

 \

 \noindent \textbf{Proof.}  The proof follows directly from Lemma 2.3 and the fact that every
linear functional is homogeneous.

\

\noindent \textbf{Lemma 2.4.}  \textit{Let $f\in Lin(M)$, and let $ g, h: M
\rightarrow \mathbb{R}$, where $M$ is a subspace of $E$,
$g(\underline 0)= h(\underline 0)= 0$, and $g(x)\le h(x)$ for every
$x\in M$. Then the following are equivalent:
\newline (i) $g(x)\le f(x)\le h(x)$ for every $x\in M$,
\newline (ii) $f(x)\le h(x+y)-g(y)$ for all $x, y\in M$.}

\

\noindent \textbf{Proof.}  First, we assume that  $g(x)\le f(x)\le
h(x)$ for every $x\in M$, and we want to show that  $f(x)\le
h(x+y)-g(y)$ for every $x, y\in M$. Since  $g(x+y)\le f(x+y)\le
h(x+y)$, and  $-h(y)\le -f(y)\le -g(y)$, it follows that
$g(x+y)-h(y)\le f(x+y)-f(y)\le h(x+y)-g(y)$. Since $f\in Lin(M)$,
then $g(x+y)-h(y)\le f(x)\le h(x+y)-g(y)$. Next, assume that
$f(x)\le h(x+y)-g(y)$ for every $x, y\in M$. Setting $y= \underline
0$ we see that $f(x)\le h(x)$. By setting $x=-y$, we obtain
$f(-y)\le -g(y)$, so $g(y)\le f(y)$. Therefore, $g(x)\le f(x)\le
h(x)$ for every $x\in M$.

\

\section{The Hahn-Banach theorem and the classical Hahn-Banach theorem}

\begin{theorem} (The Hahn-Banach theorem for sublinear functionals)

   Let $E$ be a real vector space, $M$ be a subspace of $E$ and $S \in Subl(E)$. Furthermore,
let $f_0 \in Lin(M)$ such that $ f_0(x)\leq S(x)$ for every $ x\in
M$ (or equivalently, $ -S(-x)\leq f_0(x)\leq S(x)$), then there
exists a  functional $ L\in Lin(E)$  such that $L|_M =f_0$ and $
L(x)\leq S(x)$ for every $ x\in E$ (or equivalently, $ -S(-x)\leq
f_0(x)\leq S(x)$) $[7,3]$.\end{theorem}

  The following theorem is the famous application of the Hahn-Banach
theorem and some mathematicians, like Kelly-Namioka[2], Rudin [7],
Koing [4], Simons [8,9] and others have proved it. In the following,
we introduce a simple and short proof for it.

\

\begin{theorem}(The classical Hahn-Banach theorem for sublinear functionals)
   Let $E$ be a real vector space and $S \in Subl(E)$, then there exists
a functional $ L\in Lin(E)$  such that $ L(x)\leq S(x)$ for every $
x\in E$.
\end{theorem}

\noindent \textbf{Proof.} If $S(x)=0$ for every $x\in E$, then the
result follows immediately by taking $L \equiv 0$ on $E$. Assume
that $S \not \equiv 0$ on $E$. Let $x_0 \in E$, let $M=\{\alpha
x_0:\alpha\in\mathbb{R}\}$, and define $f_0(\alpha x_0)=\alpha
S(x_0)$ for every $x\in M$. We have two cases:\newline\noindent
\textbf{Case 1}: $M =\{\underline 0\}$. Since $0 = f_0(\underline
0)= S(\underline 0)$, $f_0(x)\le S(x)$ for every $x\in M$.  By the
Hahn-Banach theorem (or Theorem 3.1), there exists a linear
functional $L$ on $E$ such that $L|_M =f_0$ and $L(x)\le S(x)$ for
every $x\in E$.\newline\noindent \textbf{Case 2}: $ M \neq
\{\underline 0\}$. Let $x = \alpha x_0 \in M$.  If $\alpha \geq 0$
then $f_0(x)=f_0(\alpha x_0)=\alpha f_0(x_0) \leq \alpha
S(x_0)=S(x)$. If $\alpha<0$ then let $\alpha=-t$, where $t>0$. Then
$f_0(x)=f_0(\alpha x_0)=\alpha S(x_0)= -t S(x_0) = -S(t x_0) =
-S(-\alpha x_0) = -S(-x) \le S(x)$ for every $x\in M$. Hence,
$f_0(x)\le S(x)$ for every $x\in M$. Thus, from the extension form
of the Hahn-Banach theorem, there exists a linear functional $L$ on
$E$ such that $L|_M=f_0$ and $L(x)\le S(x)$ for every $x\in E$.

\

   In the following Theorems, we mention a similar fact for superlinear functionals
and we call it the extension form of the Hahn-Banach theorem for
superlinear functionals and the classical Hahn-Banach theorem for
superlinear functionals [6].  The proofs of Theorems 3.3 and 3.4 can be found in [6].

\

\begin{theorem}(The extension form of Hahn-Banach theorem for superlinear functionals)

  Suppose that $f_0\in Lin(M)$, where $M$ is a subspace of $E$ such
that $P(x)\le f_0(x)$ for every $x\in M$, where $P \in Supl(E)$.
Then there exists a functional $L \in Lin(E)$ such that $L|_M= f_0$
and $P(x)\le L(x)$ for every $x\in E$.
\end{theorem}

\

\begin{theorem}(The classical Hahn-Banach theorem for superlinear functionals)

   Let $E$ be a real vector space and $P \in Supl(E)$, then there exists
a functional $ L\in Lin(E)$  such that $ P(x)\leq L(x)$ for every $
x\in E$.
\end{theorem}

\

\section{The sandwich theorem for sublinear and superlinear functionals}

   In this section we demonstrate three different approaches to proving that, if $S \in Subl(E)$, $P \in Supl(E)$
and $f_0 \in Lin(M)$, where $M$ is a subspace of $E$, and $P(x)\le
f_0(x)\le S(x)$ for each $x \in M$, then there is an extension
functional $L \in Lin(E)$ so that $P(x)\le L(x)\le S(x)$ for each $x
\in E$.  Throughout the this section, we let $T(x)=\inf\limits_{y\in
E}\{S(x+y)-P(y)\}$ for every $x \in E$.

\

\noindent \textbf{Lemma 4.1.}  \textit{Let $f_0$ be a linear functional of a subspace
$M$ of $E$ such that: $f_0(x)\le T(x)$ for every $x \in M$ \noindent
(4.1). Then the following conditions are satisfied:\newline $P(x)\le
S(x)$ for every $x\in E$  \noindent (4.2), \newline $P(x)\le
f_0(x)\le S(x)$ for every $x\in M$ \noindent (4.3).}

\

\noindent \textbf{Proof.}  From the inequality (4.1), we obtain that
$f_0(x)\le S(x+y)-P(y)$ for every $y\in E$ and $x \in M$.
Consequently, by setting $ x =\underline 0$, we see that $P(y)\le
S(y)$ for every $y\in M$. By setting $y =\underline 0$, we see that
$f_0(x)\le S(x)$ for every $x\in M$, and by setting $ y = -x$, we
obtain that $f_0(-y)\le -P(y)$ or $P(x)\le f_0(x)$ for every $x\in
M$. Hence $P(x)\le f_0(x)\le S(x)$ for every $x\in M$.

\

   The following example shows that the converse of the above lemma is
false.

\

\noindent \textbf{Example 4.1.}  \textit{Let $E=\mathbb{R}^2, P(x)=-| x_1| +
x_2, S(x)=\sqrt{x_1^2 + x_2^2}$, where $x=(x_1,x_2)$ and $M=\{(x_1,
0): x_1\in\mathbb{R}\}$ is an one dimensional subspace of $E$.
Suppose that $f_0: M\to \mathbb{R}$ is defined as $f_0(x) = x_1$ for
every $x=(x_1,0)\in M$. Then inequalities (4.2) and (4.3) are
satisfied but (4.1) is not.}

\

\textbf{Solution.} It easy to verify that $P\in Supl(E)$,
$S\in Subl(E)$ and $f_0$ is a linear functional on $M$. Let
$x=(x_1,x_2)\in E$.  Then $ S(x)=\sqrt{x_1^2+x_2^2}\ge x_2\ge -|
x_1| + x_2 = P(x)$. Also, if $x=(x_1,0)\in M$ then $P(x)=-| x_1| \le
x_1 = f_0(x)\le |x_1| \le \sqrt{x_1^2+x_2^2}= S(x)$, so $P(x)\le
f_0(x) \le S(x)$ for every  $x\in M$. Finally, suppose that
inequality (4.1) is true.  Then, for each $x\in M$, $f_0(x)\le
S(x+y)-P(y)$ for every $y\in E$. Let $x=(x_1,0)\in M$ where $x_1>0$.
Then $ x_1\le\sqrt{(x_1+y_1)^2+y_2^2}+| y_1| - y_2 $ for every
$y=(y_1,y_2)\in E$. In particular, if $y=(0,y_2)$ and $ y_2>0 $ then
this inequality becomes $ x_1 \le
\sqrt{x_1^2+y_2^2}-y_2=\frac{(\sqrt{x_1^2+y_2^2}-y_2)(\sqrt{x_1^2+y_2^2}+y_2)}
{(\sqrt{x_1^2+y_2^2}+y_2)}=\frac{x_1^2}{\sqrt{x_1^2+y_2^2}+y_2}$,
which approaches $0$ as $y_2 \rightarrow\infty $.  But this implies
that $x_1\le 0$, a contradiction.

\

\noindent \textbf{Lemma 4.2.}  \textit{If $P(x)= -S(-x)$ for every $x\in E$,
then $T(x)= S(x)$  for every  $x\in E$.}

\

\noindent \textbf{Proof.}  From the definition of $T$, we conclude
that $T(x)\le S(x+y)-P(y)$ for all $x,y\in E$. Thus $T(x)\le S(x)$.
Also, $S(x)\le S(x+y)+S(-y)= S(x+y)-P(y)$ for every $y\in E$, so
$S(x)\le T(x)$. Therefore $T(x)= S(x)$ for every $x\in E$, as
required.

\

\noindent \textbf{Corollary 4.1.} \textit{If $P(x)= -S(-x)$ for every $x\in
E$, then conditions (4.2) and (4.3) are satisfied if and only if
condition (4.1) is satisfied.}

\

\noindent \textbf{Proof.} Sufficiency follows from Lemma 4.1.
Conversely, if conditions (4.2) and (4.3) are satisfied then
condition (4.1) follows from Lemma 4.2.

\

\noindent \textbf{Lemma 4.3.}  \textit{Let $P(x)\le S(x)$ for every $x \in
E$.  Then $P(x) \le T(x) \le S(x)$ for every $x \in E$, and $T$ is a
sublinear functional. Furthermore, if $L$ is a linear functional on
$E$ such that $ P(x)\le L(x) \le S(x)$ for every $x \in E$ then
$L(x) \le T(x)$ for every $x \in E$.}

\

\noindent \textbf{Proof.} Since $P(u)\leq L(u)$, it follows that
$L(u)\le -P(-u)$. Hence, from the linearity of $L$ we see that
$L(u+v)\le S(v)-P(-u)$ for every $u,v\in E$. Setting $u+v=x $ and
$u=-y$, we obtain $L(x)\le S(x+y)-P(y)$. Taking the infimum over all
$y\in E$, we see that $L(x)\le T(x)$ for every $x \in E$. We next
prove that $T:E\to \mathbb{R}$ is sublinear. Since $P(y)\le S(y)\le
S(x+y)+S(-x)$, it follows that $-S(-x)\le S(x+y)-P(y)$. Taking the
infimum over all $y\in E,$ we see that $ T(x)\ge -S(-x)$, so
$T:E\to\mathbb{R}$. It follows that $T$ is positively homogeneous
since $S$ and $P$ are positively homogeneous.  Thus, to prove that
$T$ is sublinear it remains to show that $T$ is subadditive. Let
$u,v\in E$.  Then $ T(u+v) = \inf\limits_{y+z\in
E}\{S(u+v+y+z)-P(y+z)\}\le \inf\limits_{y\in E}\{S(u+z)-P(y)\} +
\inf\limits_{z\in E} \{S(v+z)-P(z)\} = T(u)+T(v)$. Thus $T$ is
subadditive and, consequently, sublinear. Since $T(x)\le
S(x+y)-P(y)$, it follows that $T(x)\le S(x+\underline
0)-P(\underline 0)$, so $T(x)\le S(x)$ for every $x\in E$.  Also,
$T(-y)\le S(\underline 0)-P(y)$, so $T(-y)\le -P(y)$ and $-T(y)\le
T(-y)\le -P(y)$, and hence $P(y)\le T(y)$ for every $y \in E$.
Therefore $P(x)\le T(x)\le S(x)$ for every $x \in E$, and the proof
is complete.

\

   In the following theorem we introduce the first method for proving
the statement above, which we call the sandwich theorem for
sublinear and superlinear functionals. Throughout this section, we
assume that $f_0$ is a linear functional on $M$, where $M$ is a
subspace of $E$.  As before, $S\in Subl(E) $ and $ P\in
Supl(E)$.

\

\begin{theorem}  There exists an extension linear functional $L$ on $E$ of the linear
functional $ f_0$ on $M$ such that $P(x)\le L(x)\le S(x)$ for every
$x \in E$ if and only if $f_0(x)\le T(x)$ for every $x\in M$.
\end{theorem}

\noindent \textbf{Proof.}  From Lemma 4.3 we know that $T \in
Subl(E)$ and $P(x) \le T(x) \le S(x)$ for every $x \in E$.  Also by
lemma 4.3, if $P(x) \le f_0(x) \le S(x)$ for every $x \in M$ then
$f_0(x) \leq T(x)$ for every $x \in M$ since $f_0$ is a linear
functional on $M$. Next, assume that $f_0(x) \leq T(x)$ for every $x
\in M$.  By lemma 4.1, $P(x) \leq S(x)$ for all $x \in E$ and $P(x)
\leq f_0(x) \leq S(x)$ for all $x \in M$.  By Theorem 3.3 there is
an extension linear functional $L$ on $E$ such that $f_0(x) = L(x)$
for every $x \in M$ and $P(x) \leq L(x)$ for every $x \in E$.  By
lemma 4.3 $L(x) \leq T(x) \leq S(x)$ for every $x \in E$.  Hence, it
follows that $P(x) \le L(x) \le S(x)$ for every $x \in E$.

\

   The following corollary shows that if we let $P(x)= -S(-x)$
for every $x\in E$, then we obtain the Hahn-Banach theorem.

\

\noindent \textbf{Corollary 4.2.}  \textit{If $P(x)= -S(-x)$ for every $x\in
E$ and $f_0$ is a linear functional on $M$ such that $f_0(x)\le
S(x)$ for every $x\in M$. Then there exists a linear functional $L$
on $E$ such that $L|_M= f_0$ and $L(x)\le S(x)$ for every $x\in
E$.}

\

\noindent \textbf{Proof.}  The proof follows directly from theorem
4.1 and corollary 4.1, as required.

    We now introduce the second proof of the Sandwich Theorem for
Superlinear and Sublinear Functionals.

\

\noindent \textbf{Lemma 4.4.}  \textit{Let $x_0\notin M, a=\sup\limits_{y\in
M}\{f_0(y)-S(y-x_0)\}, b=\inf\limits_{y\in M}\{-f_0(y)+S(y+x_0)\}$,
$c=\inf\limits_{y\in M}\{f_0(y)-P(y-x_0)\}$ and $d=\sup\limits_{y\in
M} \{-f_0(y)+P(y+x_0)\}$, and assume that $a\le b$ and $d\le c$.
Then the following statements are equivalent:\newline (i) $f_0(x)\le
S(x+y)-P(y)$ for every $x\in M$ and $y\in E$.\newline (ii) $a\le c$
and $d\le b$.}

\

\noindent \textbf{Proof.} Suppose that $f_0(x)\le S(x+y)-P(y)$ for
every $x\in M$ and $y\in E$. Then $f_0(w-z)\le S(w-x_0)-P(z-x_0)$
for every $w,z\in M$, where $x=w-z$ and $y=z-x_0$. Consequently,
from the linearity of $f_0$, we obtain that $f_0(w)-S(w-x_0)\le
f_0(z)-P(z-x_0)$ for every $w,z\in M$. Hence, $ \{f_0(w)-S(w-x_0):
w\in M\}$  is bounded above and has a least upper bound, $a$.
Similarly, we can see that there is a $c = \inf\limits_{z\in
M}\{f_0(z)-P(z-x_0)\}$. Since $f_0(w-z)\le S(w-x_0)-P(z-x_0)$ for
every $w,z\in M$, it follows that $a \leq c$. If we let $x=w-z$ and
$y=z+x_0$ in the inequality $f_0(x)\le S(x+y)-P(y)$, we obtain that
$-f_0(z)+P(z-x_0)\le -f_0(w)+S(w+x_0)$ for every $w,z\in M$ and
hence $ d\leq b$. Conversely, assume that $a\le c$ and $d\le b$.
Since $a\le c$, it follows that $\sup\limits_{w\in
M}\{f_0(w)-S(w-x_0)\}\le \inf\limits_{z\in M}\{f_0(z)-P(z-x_0)\}$
and hence  $f_0(w)-S(w-x_0)\le f_0(z)-P(z-x_0)$ for every $w,z\in
M$, and consequently $f_0(w-z)\le S(w-x_0)-P(z-x_0)$.  Thus,
$f_0(x)\le S(x+y)-P(y)$ for every $x=w-z\in M$ and $y=z-x_0\in E$.

\

\noindent \textbf{Remark 4.1.}  \textit{Condition (i) in the Lemma 4.4 is
equivalent to stating that $f_0(x)\le T(x)=\inf\limits_{y\in
E}\{S(x+y)-P(y)\}$ for every $x\in M$. Also, condition (ii) in Lemma
4.4 is equivalent to the condition that $ [a,b] \bigcap [d,c] \neq
\emptyset$.}

\

\noindent \textbf{Lemma 4.5.}  \textit{Let $x_0\notin M$ and let $M_1$ be
the the smallest linear subspace of $E$ containing $M$ and $x_0 $.
Let $f_0$ be a linear functional on $M$ such that $f_0(x)\le
S(x+y)-P(y)$ for every $x\in M$ and for every $y\in M_1$. Then there
exists a linear functional $L$ on $M_1$ such that $L |_M=f_0$ and
$P(x)\le L(x)\le S(x)$ for every $x\in M_1$.}

\

\noindent \textbf{Proof.}  From Lemma 4.1, we have $f_0(x)\le S(x)$
for every $x\in M$. Therefore for arbitrary $y_1,y_2\in M$ we have
$f_0(y_1+y_2)\le S(y_1+y_2)$. So, from the linearity of $f_0$ and
the sublinearity of $S$ we get $f_0(y_1)+f_0(y_2)=f_0(y_1+y_2)\le
S(y_1+y_2)\le S(y_1-x_0)+S(y_2+x_0)$. Thus, we obtain: \newline $
f_0(y_1)-S(y_1-x_0)\le - f_0(y_2)+S(y_2+x_0)$ \noindent
(4.4).\newline  By fixing $y_2$ and letting $y_1$ vary over $M$, we
see that $\{f_0(y_1)-S(y_1-x_0):y_1\in M\}$ is bounded above. Let
$a=\sup\limits_{y_1\in M}\{f_0(y_1)-S(y_1-x_0)\}$. By a similar
argument, we can find $b=\inf\limits_{y_2\in
M}\{-f_0(y_2)+S(y_2+x_0)\}$.  It is clear from (4.4) that $a\le b$.
Similarly, from Lemma 4.1, we know that $P(x)\le f_0(x)$ for every
$x\in M$. Therefore, for arbitrary $y_1,y_2\in M$ we have
$P(y_1+y_2)\le f_0(y_1+y_2)$. So, from the linearity of $f_0$ and
from the super linearity of $P$, we get $P(y_1+x_0+P(y_2-x_0)\le
P(y_1+y_2)\le f_0(y_1)+f_0(y_2)$. Thus, we have:
\newline   $-f_0(y_1)+P(y_1+x_0)\le f_0(y_2)-P(y_2-x_0)$ \noindent
(4.5).\newline By fixing $y_2$ and letting $y_1$ vary over $M$, we
see that $\{-f_0(y_1)+P(y_1+x_0):y_1\in M\}$ is bounded above. Let
$d=\sup\limits_{y_1\in M}\{-f_0(y_1)+P(y_1+x_0)\}$. By a similar
argument, there exists a $ c=\inf\limits_{y_2\in
M}\{f_0(y_2)-P(y_2-x_0)\}$, and  from inequality (4.5) we see that
$d\le c$. By Remark 4.1, we note that $[a,b]\cap [d,c]\ne
\emptyset$.  Let $\xi\in [a,b]\cap [d,c]$.  Note that, for every $
x\in M_1$, there is a unique $y\in M $ and $\alpha \in \mathbb{R}$
such that $ x = y + \alpha x_0$.  We define a real valued functional
$L$ on $M_1$ by setting $ L(x)=L(y+\alpha x_0) = f_0(y) +
\alpha\xi$.  We wish to verify that $L$ is linear such that
$L(x)=f_0(x)$ for every $x\in M$, and that $P(x)\le L(x)\le S(x)$
for every $x\in M_1$. First, we show that $L$ is linear.  Let
$z_1,z_2\in M_1$. Then there exist $y_1, y_2 \in M$ and $\alpha_1,
\alpha_2 \in \mathbb{R} $  such that $z_1=y_1+\alpha_1 x_0$ and
$z_2=y_2+\alpha_2 x_0$. Hence, for every $\alpha, \beta \in
\mathbb{R} $ we get $L(\alpha z_1 + \beta z_2) = L((\alpha y_1 +
\beta y_2) + (\alpha \alpha_1 + \beta \alpha_2) x_0) = f_0(\alpha
y_1 + \beta y_2) + ( \alpha \alpha_1 + \beta \alpha_2)\xi = \alpha (
f_0(y_1) + \alpha_1 \xi) + \beta ( f_0(y_2) + \alpha_2 \xi) = L(z_1)
+ L(z_2)$. Next, we wish to show that $L(x)=f_0(x)$ for every $x\in
M$.  If $x \in M$, then since $x = y + \alpha x_0$ for a unique
choice of $y \in M$ and $\alpha \in \mathbb{R}$, it follows that
that $x=y$ and $\alpha =0$. Thus, $L(x)=f_0(x)$ for every $x\in M$,
so $L|_M=f_0$. To show that $P(x)\le L(x)\le S(x)$ for every $x = y
+ \alpha x_0 \in M_1$, we consider three cases.
\newline {\bf Case 1.} $\alpha = 0$.  We have already shown that $L(x)=f_0(x)$,
and therefore $P(x)\le L(x)\le S(x)$, for every $ x \in M$.\newline
{\bf Case 2.}  $\alpha > 0$. Since $\xi\in [a,b]\cap [d,c]$, it
follows that $d=\sup\limits_{y\in M}\{-f_0(y)+P(y+x_0)\}\le\xi\le
b=\inf\limits_{y\in M}\{-f_0(y) + S(y+x_0)\}$.  Thus,
$-f_0(y)+P(y+x_0)\le\xi\le -f_0(y)+S(y+x_0)$ for every $y\in M$.
Replacing $y$ by $\frac{y}{\alpha}$, and multiplying the inequality
by $\alpha$, we get $P(y+\alpha x_0)\le f_0(y) +\alpha\xi \le
S(y+\alpha x_0)$ for every $y\in M$. Therefore $P(x)\le L(x)\le
S(x)$ for every $x\in M_1$, where $\alpha > 0$.\newline {\bf Case
3.} $\alpha < 0$. Let $\alpha = -t$, where $t>0$. We know that
$$a=\sup_{y\in M}\{f_0(y)-S(y-x_0)\}\le \xi\le c = \inf_{y\in
M}\{f_0(y) - P(y-x_0)\} $$  so $f_0(y)-S(y-x_0)\le\xi\le f_0(y)-
P(y-x_0)$ for every $y\in M$. Thus, $-S(y-x_0)\le - f_0(y) +\xi \le
- P(y-x_0)$ for every $y\in M$.  Replacing $y$ by $\frac{y}{t}$, and
multiplying both sides by $t$, we get $-S(y - tx_0)\le -f_0(y) +
t\xi\le -P(y - tx_0)$, so $-S(y + \alpha x_0)\le - f_0(y) - \alpha
\xi \le -P(y + \alpha x_0)$. This means that $P(y+\alpha x_0)\le
f_0(y)+\alpha\xi \le S(y+\alpha x_0)$ for every $y\in M$. Therefore,
$P(x) \le L(x) \le S(x)$ for every $x\in M_1$. This completes the
proof.

\

\begin{theorem} Let $M$ and $E_1$ be two linear subspaces of $E$ such that
$M\subseteq E_1\subseteq E$. Suppose that $f_0$ is a linear
functional on $M$ such that $f_0(x)\le S(x+y)-P(y)$ for every $x\in
M$ and $y\in E_1$. Then there exists a linear functional $L$ on
$E_1$ such that $L |_M = f_0$ and $P(x)\le L(x) \le S(x)$ for every
$x\in E_1$.
\end{theorem}

\noindent \textbf{Proof.}  The proof follows directly from Lemma 4.5
and Zorn's Lemma (See [3], p. 134 for this argument).

\

   It is worth pointing out that such an extension of a linear
functional $f_0$ on $M$ does not always exist on the whole space
$E$. If (4.1) is not valid on the whole space $E$ but is valid on
some linear subspace $E_1$ containing $M$, then by using the above
theorem we can extend $f_0$ to a linear functional $L$ on $E_1$ so
that $L |_M = f_0$ and $P(x)\le L(x)\le S(x)$ for every $x\in E_1$,
but it does not follow that we can extend $f_0$ to $E$.

\

\noindent \textbf{Example 4.2.}  \textit{Let $S(x)=\sqrt{\sum\limits_{i=1}^4
x_i^2}$ and $P(x)=-\sqrt{\sum\limits_{i=1}^3 x_i^2}+x_4$ for every
$x\in \mathbb{R}^{4}$. Let $M=\{(x_1,0,0,0):x_1\in \mathbb{R}\}$ be
a subspace of $\mathbb{R}^{4}$, and define a linear functional $f_0$
on $M$ as $f_0(x)=x_1$ for every $x=(x_1, 0,0,0)\in M$. Then we show
that there exists a linear functional $L$ on a subspace $E_1=
\mathbb{R}^3$ such that $L |_M=f_0$ and $P(x)\le L(x)\le S(x)$ for
every $x\in\mathbb{R}^3$, but there does not exist a linear
functional $L$ on $\mathbb{R}^4$ such that $L |_M=f_0$ and $P(x)\le
L(x)\le S(x)$ for every $x\in\mathbb{R}^4$.}

\

\noindent \textbf{Solution.} It is easy to verify that $S\in
Subl(\mathbb{R}^4)$ and $P\in Supl(\mathbb{R}^4)$. First, we show
that $f_0(x)\le S(x+y)-P(y)$ for every $x\in M$ and for every $y\in
M_1=\{(x_1,x_2,x_3,0):x_1,x_2,x_3\in\mathbb{R}\}$. Let
$x=(x_1,0,0,0)\in M$ and $y=(y_1,y_2,y_3,0)\in M_1$. Then
$f_0(x)=x_1$, $S(x+y)=\sqrt{(x_1+y_1)^2+y_2^2+y_3^2}$ and
$P(y)=-\sqrt{y_1^2+y_2^2+y_3^2}$. Therefore $f_0(x)= x_1 \leq |x_1|
= |x_1|-|y_1|+|y_1|\le|x_1+y_1|+|y_1|\le
\sqrt{(x_1+y_1)^2+y_2^2+y_3^2} +\sqrt{y_1^2+y_2^2+y_3^2}\le
S(x+y)-P(y)$. Hence, by Theorem 4.2, there exists a linear
functional $L_1$ on $M_1$ such that $L_1 |_M =f_0$ and $P(x)\le
L_1(x)\le S(x)$ for every $x\in M_1$. Suppose there is a linear
functional $L$ on $\mathbb{R}^4$ such that $L |_M=f_0$ and $P(x)\le
L(x)\le S(x)$ for every $x\in\mathbb{R}^4$. Then by Theorem 4.1 and
Remark 4.1, it must follow that $f_0(x)\le S(x+y)-P(y)$ for every
$x\in M$ and for every $y\in \mathbb{R}^4$. But if we let
$x=(10,0,0,0)\in M$ and $y=(0,0,0,1)\in \mathbb{R}^4$, then
$f_0(x)=10$, $S(x+y)=\sqrt{101}, P(y)=1$ and consequently
$S(x+y)-P(y)\cong 9.489 < 10 = f_0(x)$, which contradicts the
assumption that $f_0(x)\le S(x+y)-P(y)$.

\

\begin{theorem}  Let $f_0$ be a linear functional on $M$ such that $f_0(x)\le
S(x+y)-P(y)$ for every $x\in M$ and $x\in E$. Then there exists a
linear functional $L$ on $E$ such that $L|_M=f_0$ and $P(x)\le
L(x)\le S(x)$ for every $x\in E$.
\end{theorem}

\noindent \textbf{Proof.} The theorem follows directly from Theorem
4.2, setting $ E = E_1$.

\

   We now move to our third proof of the Sandwich Theorem for
Superlinear and Sublinear Functionals. Throughout the remainder of
this paper, let $ \Gamma_1 = \{L \in Lin(E):L|_M=f_0 $ and $L(x)\leq
S(x)$ for every $x\in E\}$ and  let $\Gamma_2 = \{L \in
Lin(E):L|_M=f_0 $ and $P(x)\leq L(x)$ for every $x\in E\}$.

\

\noindent \textbf{Remark 4.2.}  \textit{Let $f_0 \in Lin(M)$, $P\in
Supl(E)$, and $S\in Subl(E)$, so that $P(x)\le f_0(x)\le S(x)$ for
every $x\in M$, and $P(x)\le S(x)$ for every $x\in E$. Then by
Theorem 3.1 and see Theorem 3.3, it follows that $\Gamma_1$ and
$\Gamma_2$ are non-empty.}

\

\begin{theorem} A necessary and sufficient condition for $ \Gamma_1 \cap \Gamma_2
\neq \phi$ is that $f_0(x)\le T(x)$ for every $x\in M$, where
$T(x)=\inf\limits_{y\in E}\{S(x+y)-P(y)\}$.
\end{theorem}

\noindent \textbf{Proof.}  Suppose there is a linear functional $L
\in \Gamma_1 \cap \Gamma_2 $.  Then $L|_M=f_0$, $L(x)\le S(x)$ and
$P(x)\le L(x)$ for every $x\in E$. Since $P(u)\le L(u)$, $P(-u)\le
L(-u)=-L(u)$, and $L(u)\le -P(-u)$ for every $u\in E$. Therefore,
from the linearity of $L$, we have  $L(u+v)= L(u)+L(v)\le
S(v)-P(-u)$. This means that $L(x)\le S(x-u)-P(-u)$, so $L(x)\le
S(x+y)-P(y)$, where $ x = u+v\in E$ and $y = -u\in E$. Since
$L|_M=f_0$, it follows that $f_0(x)\le S(x+y)-P(y)$ for every $x\in
M$ and $y\in E$. Taking the infimum over all $y\in E$, $f_0(x)\le
T(x)$ for every $x\in M$. Conversely, assume that $f_0(x) \leq T(x)$
for every $x \in M$.  It is easy to verify that $ T\in Subl(E)$ (see
Lemma 4.3). By lemma 4.1, $P(x) \leq S(x)$ for all $x \in E$ and
$P(x) \leq f_0(x) \leq S(x)$ for all $x \in M$.  By Theorem 3.3
there is an extension linear functional $L$ on $E$ such that $f_0(x)
= L(x)$ for every $x \in M$ and $P(x) \leq L(x)$ for every $x \in
E$.  By lemma 4.3 $L(x) \leq T(x) \leq S(x)$ for every $x \in E$.
Hence, it follows that $P(x) \le L(x) \le S(x)$ for every $x \in E$.
Hence, $L\in \Gamma_2 \cap \Gamma_1$.

\

\begin{theorem} Let $f_0$ be a linear functional on $M$, where $M$ is a subspace of $E$
such that $P(x)\le f_0(x)\le S(x)$ for every $x\in M$, where $S \in
Subl(E)$ and $P \in Supl(E)$. Then the necessary and sufficient
condition that there exists a linear functional $L$ on $E$ such that
$L|_M=f_0$ and $P(x)\le L(x)\le S(x)$  for every $x\in E$ is that
$f_0(x)\le T(x)$ for every $x\in M$.
\end{theorem}

\noindent \textbf{Proof.} Since $f_0$ is a linear functional on $M$
and $P(x)\le f_0(x)\le S(x)$ for every $x\in M$, then $P(x)\le
f_0(x)$ and  $f_0(x)\le S(x)$, then from Remark 4.2, we conclude
that $\Gamma_1 \neq \emptyset$ and $\Gamma_2 \neq \emptyset$. Assume
that $L$ is a linear functional on $E$ such that $L|_M=f_0$ and
$P(x)\le L(x)\le S(x)$  for every $x\in E$. Then $\Gamma_1 \cap
\Gamma_2 \neq \phi$. Hence from Theorem 4.4, we have $f_0(x)\le
T(x)$ for every $x\in M$. Conversely, let $f_0(x)\le T(x)$ for every
$x\in M$. Then $T \in Subl(E)$, and by Theorem 4.4 we obtain that
$\Gamma_1 \cap \Gamma_2 \neq \phi$, and consequently there exists a
linear functional $L$ on $E$ such that $L|_M=f_0$ and $P(x)\le L(x)$
for every $x\in E$ and $L(x)\le S(x)$ for every $x\in E$.

\begin{center}

 REFERENCES

\end{center}

\

1.  A.T.Diab, N. Faried and Afaf Abou Elmatty, \textit{Extending
Linear Continuous Functionals with Preservation of Positivity},
Proc. Math. Phys. Soc. Egypt, 84, 2, p.p. 125-132, (2006).

2.  J. L. Kelley, I. Namioka, and co-authors, \textit{Linear
Topological Spaces}, D.Van Nostrand  Co., Inc.,
Princeton-Toronto-London-Melbourne(1963).

3.  A. N. Kolmogorov and S. V. Fomin, \textit{Introductory Real
Analysis}, Translated and Edited by Richard A. Silverman, Dover
Publications, Inc. New York (1970).

4.  H. Konig, \textit{Some Basic Theorems in Convex Analysis}, in
"Optimization and operations Research", edited by B. Korte,
North-Holland (1982).

5.  E. Muhamadiev and A. T. Diab, \textit{On extension of positive
linear functionals}, Internat. J. Math., Math. Sci. 23, 1 (2000).

6.  A. S. Okb-El-Bab, A. T. Diab and R. M. Boshnaa,
\textit{Generalization of the Hahn-Banach theorem}, Proceedings of
International Conference on a theme " The actual problems of
mathematics and computer science", Baku, 14-15.01.2009.

7.  W. Rudin,  \textit{Functional Analysis},  Mc Graw-Hill, New York
(1973).

8.  S. Simons, \textit{A new version of the Hahn-Banach theorem},
Arch. Math., 80, p.p. 630-646, (2003).

9.  S. Simons, \textit{From Hahn-Banach to Monotonicity}, Lecture
Notes in Mathematics, 1693, Second Edition, Springer Science +
Business Media B.V., 1998 Springer-Verlag, Berlin Heidelberg (2008).

\end{document}